\documentclass{article}

\usepackage{support}

\usepackage[utf8]{inputenc} 
\usepackage[T1]{fontenc}    
\usepackage{hyperref}       
\usepackage{url}            
\usepackage{booktabs}       
\usepackage{amsfonts}       
\usepackage{nicefrac}       
\usepackage{microtype}      
\usepackage{lipsum}
\usepackage{graphicx}
\graphicspath{{media/}}     
 \usepackage{enumitem}

\usepackage{amsmath}
\usepackage{amssymb}
\usepackage{algcompatible}

\usepackage{algpseudocode}

\usepackage{float}

\usepackage{algorithm}
  \usepackage{multicol}
\usepackage{multirow}
\usepackage{array}
    \newcolumntype{P}[1]{>{\centering\arraybackslash}p{#1}}
    \newcolumntype{M}[1]{>{\centering\arraybackslash}m{#1}}
\title{Robust unrelated parallel machine scheduling problem with interval release dates}

\author{ {\hspace{1mm}Mirosław Ławrynowicz, Jerzy Józefczyk}\\
 Department of Computer Science and Systems Engineering\\
	Wrocław University of Science and Technology\\
	Wybrzeże Wyspiańskiego 27, 50-370 Wrocław, Poland \\
	\texttt{\{miroslaw.lawrynowicz, jerzy.jozefczyk\}@pwr.edu.pl},  \\
	  \\
	
}

\hypersetup{
pdftitle={Robust unrelated parallel machine scheduling problem with interval release dates},
pdfauthor={Mirosław Ławrynowicz, Jerzy Józefczyk},
pdfkeywords={Job scheduling, Interval release dates, Unrelated parallel machines,  Makespan, Minimax regre, Worst-case analysis},
}

\floatname{algorithm}{}
\newcommand\floor[1]{\lfloor#1\rfloor}
\newcommand{\ceil}[1]{\left\lceil #1 \right\rceil}

\begin{document}
\maketitle

\begin{abstract}
This paper presents a profound analysis of the robust job scheduling problem with uncertain release dates on unrelated machines. Our model involves minimizing the worst-case makespan and interval uncertainty where each release date belongs to a well-defined interval.  Robust optimization requires scenario-based decision-making. A finite subset of feasible scenarios to determine the worst-case regret (a deviation from the optimal makespan) for a particular schedule is indicated. We formulate a mixed-integer nonlinear programming model to solve the underlying problem via three (constructive) greedy algorithms. Polynomial-time solvable cases are also discussed in detail. The algorithms solve the robust combinatorial problem using the makespan criterion and its non-deterministic counterpart. Computational testing compares both robust solutions and different decomposition strategies. Finally, the results confirm that a decomposition strategy applied to the makespan criterion is enough to create a competitive robust schedule.
\end{abstract}


\section{Introduction}
Robust decision-making has gained remarkable attention in the field of combinatorial optimization. Job scheduling problems under scenario-based uncertainty are an essential part of decision theory due to its applicability (e.g., 
\cite{atkinson2016robust}, \cite{NESBITT2021340},  \cite{peters2018robust}). Absolute robustness, robust deviation and relative robustness criteria are mainly used to hedge against parameters uncertainty \cite{kouvelis2013robust}. We study the minimax regret (robust deviation) criterion applied to the unrelated parallel machine scheduling problem $R|r_j|C_{\max\limits}$ under the assumption that each release date belongs to a well-defined intervals. Such an approach handles the problem through risk-averse because it minimizes the worst observed deviation from an optimal solution. Our model does not specify any prior knowledge (prior distribution) about an uncertain parameter. Hence, the unfavorable values determine the continuous and uniformly distributed intervals. A detailed description of alternative approaches for the representation of uncertainty can be found in \cite{hirshleifer1979analytics}, \cite{ayyub2006uncertainty},  \cite{bubnicki2013analysis}.

The deterministic unrelated parallel scheduling to minimize makespan $R|C_{\max\limits}$ is widely studied in the literature \cite{graham1979optimization}, \cite{pfund2004survey}, \cite{pinedo2012scheduling}. NP-hardness of $R|C_{\max\limits}$ implies the application of constant ratio approximation algorithms (e.g., \cite{lenstra1990approximation}, \cite{gairing2007faster}) or heuristic-based solutions: Tabu Search \cite{srivastava1998effective}, Recovering Beam Search \cite{ghirardi2005makespan}. Consequently, much research effort has been invested in designing algorithms for the scheduling problem with release dates $R|r_j|C_{\max\limits}$. In \cite{lin2013particle}, the author employs the Particle Swarm Optimization metaheuristic and formulates two functions to assess a lower bound for the makespan criterion. Also, efficient constructive algorithms based on dynamic programming and local optimization are given in \cite{mezentsev2019implementation}.

Even though a large body of literature refers to robust optimization, relatively few papers include only release dates as an imprecise parameter. The paper \cite{bachtler2020robust} addresses the profound analysis of the robust scheduling problem $1|r_j|C_{\max\limits}$ with uncertain release dates. Absolute robustness and robust deviation criteria are considered and the authors proved that both problems could be solved in polynomial time. The mathematical model of the scheduling-location (ScheLoc) problem assumes that the bounds of intervals depend on the jobs locations on a tree as shown in \cite{krumke2020robust}. The goal is to find a location of the machine on tree and a schedule to minimize the makespan value in the worst-case. The interval data results from the transfer of all jobs to a single machine across the tree with uncertain edge weights. The latter is proven to be polynomially solvable. Both \cite{bachtler2020robust} and \cite{krumke2020robust} use the concept of gamma-robustness to model uncertainty.\linebreak On the contrary, the authors of \cite{yue2018robust} investigate the robust (minimax regret) version of $1|r_j|WT_{\max\limits}$, which is NP-hard. To obtain a solution, they propose the modified Gusfield’s algorithm and constructive heuristic. Despite different criteria in \cite{bachtler2020robust} and \cite{yue2018robust}, it is showed that the number of feasible scenarios equals the number of given jobs.

Stochastic uncertainty model where the release dates are random variables with  arbitrary distributions is considered together with other stochastic parameters, e.g., processing times. \cite{liu2020parallel} presents such a complex ScheLoc problem with uncertain release dates and processing times. The two-stage formulation includes minimizing  the cost of jobs assignment to machines located within some region and the expected penalty cost of jobs earliness and tardiness. The same uncertain parameters are introduced in \cite{pinedo1983stochastic} where single machine scheduling to minimize the expected weighted sum of job completion times is given. \par

In general, a substantial number of scheduling problems with uncertain or imprecise parameters are inapproximable within any constant factor (e.g., \cite{jozefczyk2013scatter}, \cite{kasperski2012parallel}). We are forced to apply the probabilistic methods, branch and bound or constructive algorithms. Although modeling uncertainty in the input data has a practical impact on solution robustness, it appears to be computationally hard.  Then, it is worth conducting the theoretical analysis to create more reliable or exact solutions. Clearly, by indicating the polynomial-time solvable instances or using structural problem features, significant computational or decision benefits can be achieved.

The purpose of this paper is to provide an exhaustive study of the minimax regret scheduling problem $R|r_j|C_{\max\limits}$ with interval release dates. This problem has not been formally addressed yet in the literature to the best of our knowledge. Our contribution is threefold. Firstly, we prove that the worst-case regret can be determined for a finite set of feasible scenarios. Secondly, the following algorithms are developed to solve the robust problem: 
\begin{enumerate}[leftmargin=0.5cm]
  \item greedy (constructive) \textit{Partial\_Makespan algorithm} use only the makespan criterion (deterministic approach),
  \item greedy (constructive) \textit{Partial\_Regret} and \textit{Partial\_Regret\_Extended algorithms} use the minimax regret criterion.
\end{enumerate}
Our research also takes an overview of polynomial-time solvable cases. We have extended the methods known from the literature to assess the lower bound for the deterministic problem $R|r_j|C_{\max\limits}$. Thirdly, our primary finding is that the makespan criterion and decomposition strategy are enough to create a competitive solution. Moreover, the numerical evaluation explores how the different decomposition strategies affect the robust schedule quality.

The remainder of this article is as follows. The paper starts with the problem formulation in Section 2. The central part of the considerations given in Section 3 refers to the analysis of the worst-case scenario and feasible scenarios. Section 4 presents a meticulous analysis of the developed algorithms. Then in Section 5, the computational experiments are conducted, which is followed by the concluding remarks.

\section{Problem statement}
Let us consider a set $J=\{1,2,...,j,...,n\}$ of $n$ independent, non-preemptive jobs to be processed on a set\linebreak $M=\{1,2,...,i,...,m\}$ of $m$ unrelated machines. Denote by $p_{i,j}>0$ a processing time of the job $j$ to run on the machine $i$. All the processing times are stored in a matrix $p = \big[p_{i,j}\big]_{\substack{i=1,2,...,m \\ j=1,2,...,n}}$. The release date $r_{j}\ge 0$ of the job $j$ is imprecise and belongs to a closed interval $R_j=[r_j^{-},  r_j^{+}]$, $0 \le r_j^{-} < r_j^{+}$. Thus, uncertainty can be modeled through a scenario set $R=R_1 \times ... \times R_j \times ... \times R_n$ corresponding to the Cartesian product of well-known intervals. Each feasible scenario is expressed as a vector $r=[r_1,...,r_j,...,r_n]^{\text{T}}$, $r \in R$.

Let us define a binary decision $x_{i,k,j} \in \{0, 1\}$ which indicates if the job $j$ is $k$th, $k=1,2,...,n$, in a sequence deployed on the machine $i$. Consequently, a decision matrix $x = \big[x_{i,k,j}\big]_{\substack{i=1,2,...,m \\ k, j=1,2,...,n}}$ represents a schedule. The completion time of the $k$th job sequenced on $i$ can be calculated by the following recursive formula:
\begin{equation}
C_{i,k}(x,r)=\sum_{j=1}^n x_{i,k,j} \Big(p_{i,j} + \max\limits\big\{C_{i,k-1}(x,r) , r_j\big\} \Big), \ \ i=1,2,...,m, \ \ k=1,2,...,n, \ \ C_{i,0}(x,r) = 0,
\end{equation}
and $\forall_{\substack{i =1,2,...,m \\ k =1,2,...,n}}$ $C_{i,k}(x,r)=0$ if $x_{i,k,j}=0$. As a consequence, the makespan is obtained by:
\begin{equation}
C_{\max\limits}(x,r)=\max\limits_{i=1,2,...,m} C_{i,n_i}(x,r),
\end{equation}
where $n_i = \sum_{k=1}^n \sum_{j=1}^n x_{i,k,j}$ is the $i$th machine sequence length. For given $x$ and $r$, the regret is defined by the difference:
\begin{equation}
Q(x,r)=C_{\max\limits}(x,r) - C_{\max\limits}(x_r^{*}),
\end{equation}
where $x_r^{*}$ is the optimal schedule for the deterministic problem $R|r_j|C_{\max\limits}$, under the scenario $r \in R$. Therefore, \textit{the worst-case regret}, over all scenarios, takes the form:
\begin{equation}
Z(x)=\max\limits_{r \in R} Q(x,r),
\end{equation}
and \textit{the worst-case scenario} is equal to $r^*(x)= \text{arg}\max\limits\limits_{r \in R} Q(x,r)$.

To sum up, for given $J$, $M$, $R$, $p$, the non-deterministic scheduling problem consists in the determination of the optimal matrix $x^{*}$ to minimize:
\begin{equation}
Z(x^*)=\min\limits_{x} Z(x).
\end{equation}

\noindent The following constraints are imposed on the decision matrix $x$:
\begin{equation}
\sum_{i=1}^{m} \sum_{k=1}^{n} x_{i,k,j} = 1, \ \ \ j=1,2,...,n,
\end{equation}
\begin{equation}
\sum_{j=1}^{n} x_{i,k,j} \le 1, \ \ \ i=1,2,...,m, \ \ k=1,2,...,n,
\end{equation}
\begin{equation}
x_{i,k+1,j} - x_{i,k,j'} \le 0, \ \ \ i=1,2,...,m, \ \ k=1,2,...,n-1, \ \ j \ne j',
\end{equation}
\begin{equation}
x_{i,k,j} \in \big\{0,1\big\}.
\end{equation}

Constraint (6) points out that each job must be assigned exactly once. Next, each position can be occupied by at most one job (7). According to (8), a schedule on a single machine has to be represented by a sequence of consecutive non-zero entries of the matrix $x$. Such a technical constraint (8) improves the readability of $x$. Finally, the binary decision variable is given in (9).

For $\forall_{j = 1,2,...,n} \ r_{j}^{-} = r_{j}^{+}$, the deterministic counterpart $R|r_j|C_{\max\limits}$ is at least NP-hard \cite{pinedo2012scheduling}. Hence, the non-deterministic problem (5) is also at least NP-hard. The inapproximability of (5) is justified in the following property. 

\noindent\textit{\textbf{Property 1.} The minimax regret scheduling problem} (5) \textit{cannot be approximated within any constant factor.}

\textbf{Proof.} Let $x^*$ be the optimal schedule for the minimax regret version of $R|r_j|C_{\max\limits}$. Let us introduce an\linebreak $\alpha$-approximation algorithm, which returns the solution $x_s$ such that $Z(x^*)\le \alpha Z(x_s)$. It is enough to see that the assumption $\forall_{j = 1,2,...,n} \ r_{j}^{-} = r_{j}^{+}$ implies equality $Z(x^*) = Z(x_s) = 0$ because the set $R$ contains only  a single feasible scenario $r \in R$, $|R|=1$. Then, an approximation algorithm ensures $C_{\max\limits}(x_s,r) = C_{\max\limits}(x^{*})$ and an\linebreak optimal solution $x_s = x^*$ could be found for the parallel scheduling problem $R|r_j|C_{\max\limits}$. However, it cannot be done in polynomial time unless P=NP. \hfill $\blacksquare$ 

\section{Worst-case scenario analysis}
The number of feasible scenarios has a major impact on computational efficiency. Throughout this section, we clarify how to restrict the cardinality of $R$. The following theorem indicates that the set $R$ of feasible scenarios can be substantially reduced to the subset $\bar{R}=\big\{\bar{r}^1,...,\bar{r}^j,...,\bar{r}^n \big\}$, $\bar{R} \subseteq R$,        where $r=[r^{-}_1,...,r^{+}_j,...,r^{-}_n]^{\text{T}}$ is the extreme scenario. \medskip

\noindent \textit{\textbf{Theorem.} The worst-case scenario belongs to $\bar{R}$.}

\textbf{Proof.} It is sufficient to note that the processing times of jobs and only a single release date determine the value $C_{\max\limits}(x, r)$ under the scenario $r \in R$. More specifically, there must exist such a particular job $j$ assigned to the longest-working machine, i.e. $i= \text{arg}\max\limits_{l=1,2,...,m} C_{l,n_l}(x,r)$, $x_{i,k,j}=1$ for which the following property holds.\linebreak All jobs from the set $J^{(i,k)} = \{t \in J | x_{i,z,t},  z =k+1, k+2, ..., n_i \}$ succeeding the release date $r_j$ are sequenced without time gaps. In consequence $C_{i,k-1}(x,r) < r_j$, and the regret function can be calculated as:
\begin{equation}
Q(x,\bar{r}^j)= \Big( r^{+}_j + \sum_{t \in J^{(i,k)} \cup j} p_{i,t} \Big) - C_{\max\limits}(x_{\bar{r}^j}^{*}),
\end{equation}
where $\bar{r}^j = \text{arg}\max\limits_{\bar{r} \in \bar{R}} Q(x,\bar{r})$. Otherwise $r_j$ would not affect the makespan value. Note that if $\exists_{t \in J \setminus j} \ r_{t} \ne r_{t}^{-}$ or $r_{j} \ne r_{j}^{+}$ then the difference in (10) may only decrease, which leads to $ \max\limits_{\tilde{r} \in R \setminus \bar{R}}  Q(x,\tilde{r}) \le Q(x,\bar{r}^j)$. \hfill $\blacksquare$

We indicate an additional subset of scenarios to exclude from $\bar{R}$. Let us define the set $H(x)$ containing jobs with intervals covered by predecessors’ processing times:

\begin{equation}
H(x) = \Big\{ j \in J | C_{i, k-1}(x, r^{-}=\big[r^{-}_j\big]^{\text{T}}_{j=1,2,...,n}) \ge r_j^{+} \wedge x_{i,k,j} = 1, \ i=1,2,...,m, \ k=2,3,...,n \Big\}.
\end{equation} 

\noindent \textit{\textbf{Proposition. } The feasible scenario $\bar{r}^j$} can be excluded from $\bar{R}$ if the job $j$ belongs to (11).

\textbf{Proof.} Remark that any release date $r_j \in R_j$ does not change the makespan $C_{\max\limits}(x, \bar{r})$, $\bar{r} \in \bar{R}$, if $j \in H(x)$. The scenario $\tilde{r}=[r^{-}_1,...,r_j,...,r^{-}_n]^{\text{T}}$ in which $r_j \in R_j$ and $\forall_{t \in J \setminus j} \ r_t = r^{-}_t$ leads to $C_{\max\limits}(x, \bar{r}^j) = C_{\max\limits}(x, \tilde{r})$.\linebreak Clearly, $Q(x,\bar{r}^j) \le Q(x,\tilde{r})$ due to $C_{\max\limits}(x_{\bar{r}^j}^{*}) \ge C_{\max\limits}(x_{\tilde{r}}^{*})$, which concludes the proof. \hfill $\blacksquare$ \medskip 

\noindent As a consequence, we can immediately derive the theoretical bounds for (5):
\begin{equation}
0 \le Z(x) \le \max\limits_{j \in J \setminus H(x)} \Big\{ C_{\max\limits}(x, \bar{r}^j) - \big( r^{+}_j + \min\limits_{i=1,2,...,m} p_{i,j} \big) \Big\}.
\end{equation}

\section{Algorithms}
In this section, we propose three methods for solving the formulated minimax regret scheduling problem. Since the problem is in NP, we apply the problem-specific heuristics. Three algorithms are referred to as \textit{Partial\_Makespan} (PM), \textit{Partial\_Regret} (PR), and \textit{Partial\_Regret\_Extended} (PRE). The PM, PR, and PRE are deterministic, non-parametric, and use the different decomposition strategies that follow logically from Theorem. In addition, two polynomial solvable cases of (5) are discussed in detail.

\subsection{\textit{Partial\_Makespan algorithm} (PM)}
The PM algorithm creates a robust solution step-by-step starting from an empty schedule $x(1)$. Let us denote the current solution in the $u$th iteration, as $x(u) = \big[x_{i,k,j}(u)\big]_{\substack{i=1,2,...,m \\ k, j=1,2,...,n}}$, $u=1,2,...,n$. In each iteration, the decision-making process comprises two stages: subjective and greedy. The subjective stage involves a choice of a single job assignment. To evaluate the performance measure for a job, the following indicator is applied:
\begin{equation}
\mathrm{\Pi}_j(u) =    \begin{cases} \sum_{ l } avg_j , &   |\beta| > 1 \\ 
 \beta,  & \text{otherwise} \end{cases}
, \ \ \beta = \text{arg}\min\limits_{j \in J(u)} |U_j(\bar{r}^j; J(u))|,
\end{equation}
where the set $U_j\big(\bar{r}^j; J(u)\big) = \{t \in J(u) \setminus j | r^{-}_t < r^{+}_j \}$ contains jobs available before $j$ under $\bar{r}^j$,\linebreak $avg_j=m^{-1}\sum_{i=1}^m p_{i,j}$ is the averaged processing time of $j$ across all machines, and $j \in J(u)$ belongs to the set of jobs non-scheduled in iterations $1,2,...,u-1$. The idea behind (13) is to determine approximately how many jobs or processing time units may be executed before $r^{+}_j$. It exploits the problem structure because a value of (3) increases only if some subset of jobs is rescheduled to earlier positions in an optimal schedule. Hence, the PM selects the job $v = \text{arg}\min\limits_{j \in J(u)} \mathrm{\Pi}_j(u)$ or makes the arbitrary choice when there are more such jobs. During the greedy stage, an assignment of   minimizes the makespan:
\begin{equation}
(i', k') = \text{arg} \min\limits_{\substack{i=1,2,...,m \\ k=1,2,...,n}} C_{i, k}\big(x(u)=[x_{i,k,v}(u)=1], \bar{r}^{v}\big).
\end{equation}
The pair $(i', k')$ constitutes the binary decision $x_{i',k',v}(u)=1$. In practice, the PM minimizes the makespan for incomplete schedule $x(u)$. The constraint (8) prevents the overwrite of previously made decisions in (14) such that $v$ can only take the $(n_{i'}(u-1)+1)$th position. By applying a greedy strategy, the algorithm also decreases the difference:
\begin{equation}
Q\big(x(u),\bar{r}^v\big)=C_{\max\limits}\big(x(u),\bar{r}^v\big) - C_{\max\limits}\big(x(u)_{\bar{r}^v}^{*}\big).
\end{equation}
We deal with a single extreme scenario per iteration. This decomposition and greedy choice allow minimizing a value of (15) without any knowledge of the optimal schedule $x(u)_{\bar{r}^v}^{*}$. After the incrementation of $u$, the job $v$ is excluded from further possible decisions, i.e., $J(u+1) := J(u) \setminus v$, and the PM performs computations until $J(u) \ne \varnothing$.
\newpage
Pseudocode of \textit{the Partial\_Makespan algorithm} (PM) summarizes the described procedure.
\begin{algorithm} \caption{Algorithm PM}
\begin{algorithmic}[1]
\Require $J$, $M$, $R$, $p$ 
\Ensure $x_{\text{PM}}$ 
\State Set $u:=1$, $J(u):=J$, and the empty schedule $x(u) := \big[x_{i,k,j}(u) := 0\big]_{\substack{i=1,2,...,m \\ k, j=1,2,...,n}}$.
\State $\forall_{j \in J}$ generate the set $U_j(\bar{r}^j; J(u))$.
\State \textbf{while} $J(u) \ne \varnothing$
\State Calculate the job $v := \text{arg} \min\limits_{j \in J(u)} \mathrm{\Pi}(u)$.
\State Calculate the pair $(i', k') := \text{arg} \min\limits_{\substack{i=1,2,...,m \\ k=1,2,...,n}} C_{i, k}\big(x(u):=[x_{i,k,v}(u):=1], \bar{r}^{v}\big)$.
\State Set $x_{i',k',v}(u):=1$, $J(u+1) := J(u) \setminus v$, update sets in 2., and $u:=u+1$.
\State \textbf{end while}
\State $x_{\text{PM}} := x(n)$
\end{algorithmic}
\end{algorithm}

The initial assignments and generation of sets (Lines 1-2) require $\mathcal{O}(n\log n)$ time. Next, the subjective step (Line 4) is $\mathcal{O}(nm)$ because for each non-assigned job, the indicator and averaged processing times may be calculated in the worst case. Similarly, the makespans on $m$ machines are compared during the greedy step (Line 5) in $\mathcal{O}(nm)$ time. Incrementations are constant time (Line 6). Since the main loop (Lines 3-7) executes exactly $n$ times, the PM can be performed in $\mathcal{O}(n^2 m)$ time. The overall space complexity is $\mathcal{O}(n^2m + n^2 + m )$ where the schedule, subsets of jobs, and makespans require $\mathcal{O}(n^2m)$ (Line 1), $\mathcal{O}(n^2)$ (Line 2), and $\mathcal{O}(m)$ (Line 5), respectively.

\subsection{Lower bounds for the optimal makespan}
The PR and PRE use the minimax criterion to evaluate a solution. But due to  NP-hardness of $R|r_j|C_{\max\limits}$, it cannot be proved that $C_{\max\limits}(x,r) \ge C_{\max\limits}(x_r^{*})$ if $x_r^{*}$ is given by a non-exact algorithm. In \cite{lin2013particle}, the functions\linebreak  $LB(r)=\min\limits\limits_{j = 1,2,...,n} r_{j} + m^{-1} \sum_{j=1}^{n}\min\limits_{i = 1,2,...,m} p_{i,j}$, $LB_1(r)=\max\limits_{j = 1,2,...,n} \big( r_{j} + \min\limits_{i = 1,2,...,m} p_{i,j} \big)$ are given to assess the lower bound for $R|r_j|C_{\max\limits}$. Our research extends the bounds formulated in \cite{lin2013particle}. 

Let us introduce an instance in which the release dates are densely distributed within a single time interval. Namely, a relatively small number of jobs (greatly exceeding $m$) is available near the same time. However, the bound $LB_1(r)$ may be unachievable for some instances due to averaging operator. For example, if there exists such a release date that $LB_1(r) \le r_j$. To handle this issue, we introduce a set $E_j(r;J) = \{ t \in J | r_t \ge r_j  \}$ of jobs unavailable before $j$ and calculate:
\begin{equation}
W_j\big( E_j(r;J), p \big)= \min\limits_{j \in E_j(r;J)} r_j + m^{-1} \sum_{j \in E_j(r;J)} \min\limits_{i = 1,2,...,m} p_{i,j},
\end{equation}
The value of (16) for $j = \text{arg}\min\limits_{l = 1,2,...,n} r_{l}$ equals $LB(r)$. It is clear that the set $E_t(r;J)$, $t\ne j$, excludes some subset from $J$ which may decrease a mean processing time over $m$ machines in (16). Then, the lower bound\linebreak $LB_2(r)=\max\limits_{j=1,2,...,n} \big\{ W_j\big( E_j(r;J), p\big) \big\}$ ensures that $LB(r) \le LB_2(r)$. 

The next approach uses the assumption that the elements stored in $E_j(r;J)$, $\gamma_j = \big|E_j(r;J)\big|$, can be partitioned into  $\lambda_j= \ceil{m^{-1}\gamma_j}$ subsets of $m$ jobs simultaneously scheduled over $m$ machines. A value of $m^{-1}\gamma_j$ is rounded up because there may exist a single subset with cardinality less than $m$. Consequently, the bound is provided by adding $\lambda_j$ times the shortest processing time $p_w = \min\limits_{j \in E_j(r;J)} \min\limits_{i = 1,2,...,m} p_{i,j}$ to the earliest release date:
\begin{equation}
\tilde{W}_j\big( E_j(r;J), p_w \big)= \min\limits_{j \in E_j(r;J)} r_j + \sum_{w=1}^{\lambda_j} p_w,
\end{equation}
and, analogously, $ LB_3(r)=\max\limits_{j=1,2,...,n}\big\{ \tilde{W}_j\big( E_j(r;J), p_w \big) \big\}$.

Proposed methods are general enough to deal with unequally distributed release dates. The complementariness of selected functions results in a tight bound $ LB(\bar{r}) = \max\limits\big\{LB_1(\bar{r}),LB_2(\bar{r}),LB_3(\bar{r}) \big\}$. Lastly, the application of lower bounds to (4) defines \textit{the relaxed worst-case regret}:
\begin{equation}
\tilde{Z}(x)= \max\limits_{\bar{r} \in \bar{R}} \big\{ C_{\max\limits}(x, \bar{r}) - LB(\bar{r})\big\}.
\end{equation}
\subsection{\textit{Partial\_Regret algorithm} (PR)}
The PR algorithm adopts the greedy strategy based on a robust criterion. The developed method implements a fundamentally different optimization methodology than the PM. Namely, the decisions are set taking into account the decomposition of \textit{the relaxed worst-case regret} (18) instead of the two-step evaluation in (13)-(14). 

In the $u$th iteration, the PR search for the decision $x_{i', k', v'}(u)=1$ such that:
\begin{equation}
(i', k', v') = \text{arg} \min\limits_{\substack{i=1,2,...,m \\ k=1,2,...,n \\ v \in J(u)} } \big\{ C_{i, k}\big(x(u)=[x_{i,k,v}(u)=1], \bar{r}^{v}\big) - LB(\bar{r}^v)\big\}.
\end{equation}
\noindent Note that the difference in (19) uses the incomplete schedule $x(u)$ and considers only a single extreme scenario for each non-assigned job. Hence, the PR decomposes (18) and operates on partial regrets. The greedy approach comes down to iteratively solve $n$ subproblems given in (19). Let us introduce the set $T$ containing jobs of equivalent solutions in (19). If $|T|>1$, the following indicator is applied:
\begin{equation}
\varphi_j(x, u) = \max\limits\big\{ C_{i', k' - 1}\big(x(u-1), \bar{r}^{-}\big)- r_{j}^{+}, 0 \big\}.
\end{equation}

\noindent It expresses the time gap between the last job in a sequence and the current job $j$, under a particular scenario. 
Clearly, the job $\text{arg}\max\limits\limits_{j \in T} \varphi_j(x, u)$ constitutes a fairly tight schedule. The arbitrary choice is made when more jobs is indicated. In order to reduce the computational effort, the lower bounds $LB(\bar{r})$, $\bar{r} \in \bar{R}$, are calculated at preprocessing phase. The PR terminates when a schedule is complete.

Pseudocode of \textit{the Partial\_Regret algorithm} (PR) summarizes the described procedure.

\begin{algorithm} \caption{Algorithm PR}
\begin{algorithmic}[1]
\Require $J$, $M$, $R$, $p$ 
\Ensure $x_{\text{PR}}$ 
\State Set $u:=1$, $J(u):=J$, and the empty schedule $x(u) := \big[x_{i,k,j}(u) := 0\big]_{\substack{i=1,2,...,m \\ k, j=1,2,...,n}}$.
\State $\forall_{j \in J}$ generate the lower bound  $LB(\bar{r}^j)$.
\State \textbf{while} $J(u) \ne \varnothing$
\State Calculate $(i', k', v') := \text{arg} \min\limits\limits_{\substack{i=1,2,...,m \\ k=1,2,...,n \\ j \in J(u)} } \big\{ C_{i, k}\big(x(u):=[x_{i,k,j}(u)=1], \bar{r}^{j}\big) - LB(\bar{r}^j)\big\}$.
\State \textbf{if} $|T|>1$ \textbf{then} $v' := \text{arg}\max\limits\limits_{j \in T} \varphi_j(x, u)$ \textbf{end if}
\State Set $x_{i',k',v'}(u):=1$, $J(u+1) := J(u) \setminus v'$ and $u:=u+1$.
\State \textbf{end while}
\State $x_{\text{PR}} := x(n)$
\end{algorithmic}
\end{algorithm}
The initial assignments and generation of sets (Lines 1-2) require $\mathcal{O}(n^2)$ time. Next, the greedy decision (Line 4) is $\mathcal{O}(nm)$ (analogously to Line 5 of the PM) due to the memorized lower bounds and makespans. Since the subset $|T|=J$ in the worst case (line 5), the procedure takes $\mathcal{O}(n)$ time. Incrementations are constant time (Line 6).  Finally, the PR time complexity is $\mathcal{O}(n^2m)$ The overall space complexity is $\mathcal{O}(n^2m + n + m )$ where the schedule, makespans, and lower bounds require $\mathcal{O}(n^2m)$ (Line 1), $\mathcal{O}(m)$ (Line 2), and $\mathcal{O}(n)$ (Line 4), respectively.
\subsection{\textit{Partial\_Regret\_Extended algorithm} (PRE)}
Both the PM and PR involve each subsequent decision for a particular extreme scenario. In (14) and (19), only a single extreme scenario is considered for a non-assigned job. The PRE use each extreme scenario $\bar{r}^t$ to make a single decision where $t$ belongs to the set $t \in J \setminus \{J(u) \setminus v \}$ of assigned jobs and $v \in J(u)$. Then, the greedy choice comprises two nested problems:
\begin{equation}
(i', k', v') = \text{arg} \min\limits_{\substack{i=1,2,...,m \\ k=1,2,...,n \\ v \in J(u)} } \Big\{ \max\limits_{t \in J \setminus \{J(u) \setminus v \}} \big\{ C_{i, k}\big(x(u)=[x_{i,k,t}(u)=1], \bar{r}^{t}\big) - LB(\bar{r}^t) \big\} \Big\}.
\end{equation}
\noindent Unlike the PR and (19), (21) calculates the relaxed worst-case regret for incomplete schedule and given bounds.\linebreak It enables the verification of more greedy decisions than the PR at the cost of computational efficiency. The application of each extreme scenario of  $\bar{r}^t$, $t \in J \setminus J(u) $, in (21) forces the modification of (20). The averaged time gap for the job $j$:
\begin{equation}
\bar{\varphi}_j(x, u) = \big|J \setminus J(u)\big|^{-1} \sum_{v \in J \setminus J(u) } \max\limits\big\{ C_{i', k' - 1}\big(x(u-1), \bar{r}^{v}\big)- r_{j}^{+}, 0 \big\}
\end{equation}
indicates the best equivalent solution of (21). 

The PRE pseudocode replaces Lines 4 and 5 of the PR by (21) and (22). Other lines remain unchanged. The nested problems in (21) increase algorithm time complexity in comparison to the PM and PR. The problem (21) takes $\mathcal{O}(n^2m)$ time due to the nested maximization problem. For a single job, the indicator (22) is $\mathcal{O}(n)$ in the worst case $J(u)=\varnothing$. Then, the PRE time complexity is $\mathcal{O}(n^3m)$ time. Since the PRE only implements the modified equations (19) and (20), the space complexity is $\mathcal{O}(n^2m + n + m )$ (analogously to the PR).
\subsection{Polynomial-time solvable cases}
The hardness of the robust problem (5) is highly correlated with the interval bounds. We formulate two conditions that simplify the problem and prove that the PM yields an optimal solution for given cases. At first, let us assume the disjoint intervals:
\begin{equation}
\forall_{ \{j, t\} \subseteq J} \ R_j \cap R_t= \varnothing .
\end{equation}

\noindent \textit{\textbf{Property 2.} The schedule $x_{\text{PM}}$ is optimal if }(23) \textit{holds.}

\textbf{Proof.} Note that the order of release dates  $r_{j_1}<...<r_{j_b}<...<r_{j_n}$, $r_{j_b} \in R_{j_b}$, ensures $|U_{j_1}(\bar{r}^{j_1}; J(1))|=...=|U_{j_b}(\bar{r}^{j_b}; J(u))|=...=|U_{j_n}(\bar{r}^{j_n}; J(u))| = 0$. Hence, both makespans $C_{\max\limits}(x_{\text{PM}},\bar{r}^{j_b}) = C_{\max\limits}(x_{\bar{r}^{j_b}}^{*})$, $b=1,2,...,n$, have the lowest possible values under any scenario and, in effect, $Z(x_{\text{PM}})=0$. \hfill $\blacksquare$ \medskip

\noindent Second condition:
\begin{equation}
\exists_{ j \in J} \ \max\limits_{k \in J \setminus j} r_{k}^{+} + \sum_{k \in J \setminus j} \max\limits_{i=1,2,...,m} p_{i,k} \le r_{j}^{+} ,
\end{equation}
refers to an instance in which the optimality of the solution depends only on the position of $j$.

\noindent \textit{\textbf{Property 3.} The schedule $x_{\text{PM}}$ is optimal if }(24) \textit{holds.}

\textbf{Proof.} The inequality (24) leads to $C_{\max\limits}(x,r) \le r_{j}^{-}$ irrespective of either any scenario or schedule of $J \setminus j$. The job $j$ is always placed after any element in $J \setminus j$ according to $x_{\text{PM}}$ due to $|U_j(\bar{r}^j; J(u))| = 0$. Then, \textit{the worst-case regret}, under any scenario, takes the form:\begin{equation}
Z(x_{\text{PM}}) = r_j + p_{i,j} - (r_j + p_{i',j}) = p_{i,j} - p_{i',j}, \ \ r_j \in R_j, \ \ \{i, i'\}\subseteq M,
\end{equation}
 and the greedy choice in (14) leads to $i=i'$. \hfill $\blacksquare$ 

\section{Computational results}
This section presents a series of numerical experiments to compare the developed algorithms. All algorithms are implemented in Python 3, and the computational experiments are performed on an Apple M1 CPU with 16 GB of RAM. Before discussing our research, we will focus on dataset creation. Practical methods for generating hard instances refer mainly to interval processing times \cite{cwik2018heuristic}, \cite{allahverdi2014single}, \cite{sotskov2009minimizing}. Only \cite{yue2018robust} describes how to create a dataset for the robust problem with interval release dates. However, this approach is prone to develop accessible instances (\textit{\textbf{Properties 2-3}}) for our problem.

Some instances require solving $R|r_j|C_{\max\limits}$ instead of the robust counterpart. For example, let us introduce the subset $\tilde{J}= \{j_1,j_2,...,j_n\}$, $\tilde{J} \subset J$, of jobs satisfying (23), $|M|=m$, where the order $r^{-}_{j_1} < r^{-}_{j_2} <...< r^{-}_{j_m}<r^{-}_{j_{m+1}} \le ... < r^{-}_{j_n}$ is preserved. Since any scenario does not change the order of $\tilde{J}$ assigned according to \textit{\textbf{Property 2}}, the remaining jobs $J \setminus \tilde{J}$ are scheduled later on each machine. By assuming the condition $\forall_{j \in \tilde{J} \land l \in J \setminus \tilde{J}} \ \big( r_{j}^{-} + \min\limits_{i = 1,2,...,m} p_{i,j} \ge r_{l}^{+} \big)$, we are forced to schedule irrespective of any scenario because the processing times of jobs in $\tilde{J}$ cover all intervals $R_l$, $l \in J \setminus \tilde{J}$. To handle the above-mentioned observations, we propose two complementary datasets $DS_1$, $DS_2$ where the release date intervals are densely and sparsely distributed, respectively. More specifically, all jobs' processing times and release dates are randomly drawn from the discrete uniform distribution. Each interval is constrained such that $R_j= [ r^{-}_{j}, r^{+}_{j}] = [r^{-}_{j}, r^{-}_{j}+ avg_j * offset_j]$, $r^{-}_{j}<r^{+}_{j}$, $offset_j \sim U(0.2, 5)$, where the term $r^{-}_{j}+ avg_j * offset_j$ prevents from unreasonable long intervals. To avoid (23) and (24), we divide the timeline into $w$ consecutive and disjoint time segments and generate at least   $\floor{n/w}$ release date intervals within each segment. Clearly, a value of $w$ defines the intervals density. Our study includes the instances denoted by triple $\langle m,n,D \rangle	$ where $m \in \{2,3,...,20\}$, $n \in \{50, 100,...,500 \}$ and $D \in \{DS_1, DS_2\}$. Both datasets are prepared as follows: 
\begin{enumerate}[leftmargin=0.5cm]
  \item $DS_1$: $R_j= [ r^{-}_{j}, r^{+}_{j}] \subseteq [0, 150]$, $p_{i,j} \in [5, 50]$, $w=10$,
  \item $DS_2$: $R_j= [ r^{-}_{j}, r^{+}_{j}] \subseteq [0, 300]$, $p_{i,j} \in [5, 50]$, $w=5$.
\end{enumerate}
At first, we generate two prominent instances $I_{20, 500, k} = \langle m=20, n=500, DS_k \rangle$, $k \in \{1,2\}$. Next, the remaining instances are created by uniformly removing data from each time segment of $I_{20, 500, k}$. In consequence, all cases share same data.

The first set of experiments is carried out to evaluate the quality of schedules (Tables 1-4). The solutions obtained by the PM, PR, and PRE are denoted by $x_{\text{PM}}$, $x_{\text{PR}}$, $x_{\text{PRE}}$, respectively. In each table, the best objective value yielded\linebreak by the PM is printed in bold type.
{\renewcommand{\arraystretch}{1.4}
\begin{table}[h]
\centering
\caption{Dependence of \textit{the relaxed worst-case regret} on $n$ and $m$ ($DS_1$)}
\begin{tabular}{|c|c|c|c|c|c|c|c|c|c|} 
\hline
\multirow{2}{*}{$n$} & \multicolumn{3}{c|}{$m$=5} & \multicolumn{3}{c|}{$m$=10} & \multicolumn{3}{c|}{$m$=15}  \\ 
\cline{2-10}
                   & $\tilde{Z}(x_{\text{PM}})$ & $\tilde{Z}(x_{\text{PR}})$ & $\tilde{Z}(x_{\text{PRE}})$              & $\tilde{Z}(x_{\text{PM}})$ & $\tilde{Z}(x_{\text{PR}})$ & $\tilde{Z}(x_{\text{PRE}})$              & $\tilde{Z}(x_{\text{PM}})$ & $\tilde{Z}(x_{\text{PR}})$ & $\tilde{Z}(x_{\text{PRE}})$               \\ 
\hline
50                 & \textbf{40.0}    & 95    & 57.0       & \textbf{0}     & 43.0    & 0         & \textbf{0}    & 0     & 0           \\ 
\hline
100                & 97.4  & 130.4 & 67.4     & 23.0    & 33.0    & 18.0        & 11.0   & 27.0    & 7.0           \\ 
\hline
150                & 124.2 & 108.2 & 48.2     & 32.0    & 76.0    & 45.0        & 23.0   & 42.0    & 16.0          \\ 
\hline
200                & 145.4 & 115.4 & 54.5     & 75.5  & 116.5 & 58.0        & 17.0   & 41.0    & 16.0          \\ 
\hline
250                & 147.6 & 115.6 & 51.6     & 82.2  & 113.2 & 45.2      & 34.0   & 62.0    & 26.0          \\ 
\hline
300                & 199.0   & 104.0   & 45.0       & 88.1  & 121.1 & 52.1      & 57.0   & 91.0    & 46.0          \\ 
\hline
350                & 245.4 & 110.4 & 45.4     & 94.8  & 110.8 & 43.8      & 73.0   & 97.0    & 55.0         \\ 
\hline
400                & 248.6 & 105.6 & 44.1     & 106.2 & 113.2 & 48.2      & 73.2 & 103.2 & 44.1        \\ 
\hline
450                & 253.2 & 118.2 & 49.9     & 123.0   & 118.0   & 41.0        & 79.5 & 101.5 & 48.0          \\ 
\hline
500                & 313.0   & 97.0    & 42.4     & 118.2 & 120.2 & 42.0        & 81.5 & 97.5  & 38.5        \\
\hline
\end{tabular}
\end{table}
}

{\renewcommand{\arraystretch}{1.4}
\begin{table}[h]
\centering
\caption{Dependence of \textit{the relaxed worst-case regret} on $n$ and $m$ ($DS_1$)}
\begin{tabular}{|c|c|c|c|c|c|c|c|c|c|} 
\hline
\multirow{2}{*}{$m$} & \multicolumn{3}{c|}{$n$=150} & \multicolumn{3}{c|}{$n$=300} & \multicolumn{3}{c|}{$n$=450}  \\ 
\cline{2-10}
                    & $\tilde{Z}(x_{\text{PM}})$ & $\tilde{Z}(x_{\text{PR}})$ & $\tilde{Z}(x_{\text{PRE}})$              & $\tilde{Z}(x_{\text{PM}})$ & $\tilde{Z}(x_{\text{PR}})$ & $\tilde{Z}(x_{\text{PRE}})$              & $\tilde{Z}(x_{\text{PM}})$ & $\tilde{Z}(x_{\text{PR}})$ & $\tilde{Z}(x_{\text{PRE}})$ \\ 
\hline
2                  & 208.0   & 106.0   & 46.0         & 315.0   & 81.0    & 24.0         & 570.5 & 104.5 & 29.5        \\ 
\hline
4                  & 126.3 & 102.3 & 43.3       & 248.0   & 87.0    & 45.0         & 384.3 & 100.3 & 32.3        \\ 
\hline
6                  & 95.3  & 108.3 & 60.3       & 144.7 & 108.7 & 44.7       & 223.7 & 104.7 & 41.7        \\ 
\hline
8                  & 69.3  & 106.3 & 60.3       & 125.9 & 108.9 & 38.9       & 164.4 & 100.4 & 45.4        \\ 
\hline
10                 & \textbf{32.0}    & 76.0    & 45.0         & 88.1  & 121.1 & 52.1       & 123.0   & 118.0   & 41.0          \\ 
\hline
12                 & \textbf{19.0}    & 59.0    & 20.0         & 76.3  & 112.3 & 51.3       & 91.1  & 121.1 & 47.1        \\ 
\hline
14                 & 23.0    & 33.0    & 16.0         & 69.0    & 107.0   & 50.1        & 89.6  & 107.6 & 38.6        \\ 
\hline
16                 & 23.1    & 31.0    & 14.4         & 45.0    & 77.1    & 38.0         & 68.8  & 94.8  & 29.9        \\ 
\hline
18                 & 18.0    & 20.0    & 9.0          & \textbf{28.0}    & 53.0   & 28.0        & 75.6  & 91.6  & 36.1        \\ 
\hline
20                 & 18.0   & 18.0    & 9.0          & 31.0    & 41.0    & 26.0        & 53.0    & 92.0    & 45.0          \\
\hline
\end{tabular}
\end{table}
}
\newpage Analyzing the results in Tables 1-4, we see the conspicuous quality of $x_{\text{PRE}}$. A major reason the PRE outperforms other approaches is that it reconsiders the decisions for each position separately. In consequence, it makes considerably more iterations for each subproblem. However, there exists a subset of instances where the makespan criterion is enough to create a competitive solution. In Tables 1-2, the PM  performs well on instances where the cardinalities $U_j\big(\bar{r}^j; J(u)\big)$, $j\in J(u)$, are unequal and the value $r^{+}_{j} + p_{i,j}$, $j\in J(u)$, do not cover the significant number (greater than three) of consecutive intervals. Unfortunately, an increasing number (density) of jobs worsens the quality of $x_{\text{PM}}$ as shown in Table 1. This is because, after relatively few consecutive decisions, the processing times of scheduled jobs, under any scenario, cover the intervals of remaining jobs. Then, the PM  solves the deterministic problem because the release dates of non-scheduled jobs do not change the makespan. Please note that the replacement of (13) and (14) by the simplified regret evaluation in (19) improves the robust solution (compare $x_{\text{PM}}$ in Table 1 for m=5). Consequently, since the density of jobs within each time segment during the decision-making procedure decreases with the number of machines, the schedule $x_{\text{PM}}$ ensures the lowest value of (18) (Tables 1-2).

Based on the previous observations, we provide the same experiments for $DS_2$ (Tables 3-4). We also observe that the PM  schedules effectively for instances where the PR and PRE have to enhance their decisions by (20) and (22) due to many equivalent candidate solutions in subsequent iterations. Undoubtedly, the indicator (13) accurately determines a job order to assign so that a schedule can be robust for a wide range of instances in Tables 3-4. Analogously to Tables 1-2, it is empirically confirmed that the quality of $x_{\text{PM}}$ deviates significantly from $x_{\text{PR}}$ and $x_{\text{PRE}}$ when the density of jobs increases.
{\renewcommand{\arraystretch}{1.4}
\begin{table}[h]
\centering
\caption{Dependence of \textit{the relaxed worst-case regret} on $n$ and $m$ ($DS_2$)}
\begin{tabular}{|c|c|c|c|c|c|c|c|c|c|} 
\hline
\multirow{2}{*}{$n$} & \multicolumn{3}{c|}{$m$=5} & \multicolumn{3}{c|}{$m$=10} & \multicolumn{3}{c|}{$m$=15}  \\ 
\cline{2-10}
                    & $\tilde{Z}(x_{\text{PM}})$ & $\tilde{Z}(x_{\text{PR}})$ & $\tilde{Z}(x_{\text{PRE}})$              & $\tilde{Z}(x_{\text{PM}})$ & $\tilde{Z}(x_{\text{PR}})$ & $\tilde{Z}(x_{\text{PRE}})$              & $\tilde{Z}(x_{\text{PM}})$ & $\tilde{Z}(x_{\text{PR}})$ & $\tilde{Z}(x_{\text{PRE}})$  \\ 
\hline
50                 & 17.0    & 17.0    & 16.0       & \textbf{0}     & 0     & 0         & 1  & 0  & 0                \\ 
\hline
100                & 43.0    & 81.0    & 63.0       & 18.0    & 26.0    & 16.0        & 18.0 & 26.0 & 9.0                \\ 
\hline
150                & 114.2 & 139.2 & 104.2    & \textbf{18.0}    & 26.0    & 26.0        & 20.0 & 26.0 & 9.0                \\ 
\hline
200                & 112.4 & 142.4 & 79.4     & \textbf{20.0}    & 43.1    & 24.0        & 20.0 & 23.0 & 13.0               \\ 
\hline
250                & 146.6 & 140.6 & 78.6     & \textbf{24.0}    & 57.0    & 28.0        & 25.0 & 29.0 & 13.0               \\ 
\hline
300                & 176.0   & 136.0   & 72.0       & 38.0    & 100.0   & 36.0        & 25.0 & 31.1 & 24.0               \\ 
\hline
350                & 217.4 & 135.4 & 68.4     & 89.0    & 130.0   & 62.0        & \textbf{25.0} & 27.4 & 26.0               \\ 
\hline
400                & 238.6 & 132.6 & 66.8     & 113.2 & 148.2 & 74.2      & \textbf{25.0} & 32.0 & 26.0               \\ 
\hline
450                & 252.2 & 137.2 & 69.9     & 124.0   & 148.0   & 77.0        & 32.0 & 60.0 & 26.0               \\ 
\hline
500                & 322.0   & 137.0   & 68.1     & 135.2 & 154.2 & 79.1      & \textbf{27.0} & 71.0 & 30.0               \\
\hline
\end{tabular}
\end{table}
}
\vspace*{-.3cm}

{\renewcommand{\arraystretch}{1.4}
\begin{table}[h]
\centering
\caption{Dependence of \textit{the relaxed worst-case regret} on $n$ and $m$ ($DS_2$)}
\begin{tabular}{|c|c|c|c|c|c|c|c|c|c|} 
\hline
\multirow{2}{*}{$m$} & \multicolumn{3}{c|}{$n$=150} & \multicolumn{3}{c|}{$n$=300} & \multicolumn{3}{c|}{$n$=450}  \\ 
\cline{2-10}
                  & $\tilde{Z}(x_{\text{PM}})$ & $\tilde{Z}(x_{\text{PR}})$ & $\tilde{Z}(x_{\text{PRE}})$              & $\tilde{Z}(x_{\text{PM}})$ & $\tilde{Z}(x_{\text{PR}})$ & $\tilde{Z}(x_{\text{PRE}})$              & $\tilde{Z}(x_{\text{PM}})$ & $\tilde{Z}(x_{\text{PR}})$ & $\tilde{Z}(x_{\text{PRE}})$          \\ 
\hline
2                  & 165.0   & 104.0   & 55.0         & 340.0   & 104.0   & 56.0         & 618.5 & 110.5 & 37.5        \\ 
\hline
4                  & 189.3 & 155.3 & 102.3      & 266.0   & 154.0   & 59.0         & 337.3 & 151.3 & 56.3        \\ 
\hline
6                  & \textbf{35.0}    & 117.0   & 62.0         & 148.7 & 155.7 & 70.7       & 196.7 & 153.7 & 58.6        \\ 
\hline
8                  & \textbf{23.0}    & 41.0    & 26.0         & 130.9 & 152.9 & 89.9       & 175.4 & 151.4 & 64.4        \\ 
\hline
10                 & \textbf{18.0}    & 26.0    & 26.0         & 38.0    & 100.0   & 36.0         & 124.0   & 148.0   & 77.0          \\ 
\hline
12                 & \textbf{20.0}    & 26.0    & 26.0         & \textbf{25.0}    & 33.0    & 27.0         & 85.0    & 131.0   & 39.0          \\ 
\hline
14                 & \textbf{20.0}    & 26.0    & 20.0         & 25.0    & 31.0    & 17.0         & 29.0    & 65.0    & 19.0          \\ 
\hline
16                 & 20.0    & 26.0    & 9.0          & \textbf{25.0}    & 29.0    & 27.0         & 24.0    & 29.0    & 15.0          \\ 
\hline
18                 & 5.0     & 10.0    & 2.0          & \textbf{12.0}    & 15.0    & 13.0         & 30.0    & 24.2    & 13.0        \\ 
\hline
20                 & 5.0    & 5.0     & 1.0          & \textbf{10.0}    & 10.0    & 10.0         & 17.0    & 21.0    & 13.0          \\
\hline
\end{tabular}
\end{table}
}

In Table 5, we present the running time comparison for $DS_1$. The results for $DS_2$ are similar to those shown in Table 5. The running times of the PM, PR, and PRE are denoted by $t_{\text{PM}}$, $t_{\text{PR}}$, $t_{\text{PRE}}$, respectively. The major drawback of the PRE is the considerably long time of computation. It results from the formulation of nested problems in (21). On the other hand, the PM is extremely fast, even for large instances. The difference in the running times between the PM and PR comes from calculating lower bounds and evaluation of (19). 

{\renewcommand{\arraystretch}{1.4}
\begin{table}[h]
\centering
\caption{Dependence of the running time on $n$ and $m$ ($DS_1$)}
\begin{tabular}{|c|c|c|c|c|c|c|c|c|c|} 
\hline
\multirow{2}{*}{$n$} & \multicolumn{3}{c|}{$m$=5} & \multicolumn{3}{c|}{$m$=10} & \multicolumn{3}{c|}{$m$=15}  \\ 
\cline{2-10}
                   & $t_{\text{PM}}$ [ms]    & $t_{\text{PR}}$ [ms]      & $t_{\text{PRE}}$ [s]             &    $t_{\text{PM}}$ [ms]    & $t_{\text{PR}}$ [ms]      & $t_{\text{PRE}}$ [s]              &    $t_{\text{PM}}$ [ms]    & $t_{\text{PR}}$ [ms]      & $t_{\text{PRE}}$ [s]               \\ 
\hline
50                 & 9   & 254  & 3           & 9   & 254  & 6            & 9   & 262  & 7             \\ 
\hline
100                & 32  & 753  & 21          & 32  & 799  & 38           & 28  & 893  & 51            \\ 
\hline
150                & 78  & 1522 & 65          & 78  & 1572 & 116          & 79  & 1592 & 185           \\ 
\hline
200                & 121 & 1710 & 147         & 121 & 1801 & 264          & 143 & 2213 & 445           \\ 
\hline
250                & 242 & 3093 & 303         & 242 & 3063 & 599          & 243 & 3194 & 957           \\ 
\hline
300                & 333 & 3901 & 521         & 333 & 3975 & 903          & 329 & 4007 & 1413          \\ 
\hline
350                & 386 & 4876 & 762         & 384 & 4899 & 1519         & 434 & 5020 & 2029          \\ 
\hline
400                & 518 & 5600 & 982         & 511 & 5864 & 2345         & 502 & 6387 & 3308          \\ 
\hline
450                & 725 & 6647 & 1200        & 725 & 6596 & 3151         & 764 & 7359 & 4880          \\ 
\hline
500                & 910 & 9179 & 1401        & 916 & 9399 & 4211         & 921 & 9835 & 5891          \\
\hline
\end{tabular}
\end{table}

}
Note that both the PR and PRE use a priori given $LB(\bar{r}^j)$, $j=1,2,...,n$, for the complete set of jobs $j \in J$ and extreme scenarios $\bar{r}^j \in \bar{R}$. Let us modify (19) and (21) such that the lower bounds depend only on the jobs scheduled according to the incumbent solution $x(u)$, $u=1,2,...,n$. Namely,  $j \in J\setminus J(u)$ and an additional job to be scheduled constitute the lower bound in each iteration. Then, the algorithms take into account only a subset of jobs, which constraints the planning horizon. Unfortunately, the short-sighted decision-making strategy does not improve the solutions compared to Table 1 and Table 3 as shown in Table 6. The deviations from the values in Tables 1-4 are also preserved for $m \ne 5$.

{\renewcommand{\arraystretch}{1.4}
\begin{table}[h]
\centering
\caption{Dependence of the relaxed worst-case regret on $n$ and $m=5$}
\begin{tabular}{|c|c|c|c|c|} 
\hline
\multirow{2}{*}{$n$} & \multicolumn{2}{c|}{$DS_1$} & \multicolumn{2}{c|}{$DS_2$}  \\ 
\cline{2-5}
                   &$\tilde{Z}(x_{\text{PR}})$ & $\tilde{Z}(x_{\text{PRE}})$               &       $\tilde{Z}(x_{\text{PR}})$ & $\tilde{Z}(x_{\text{PRE}})$                \\ 
\hline
50                 & 126.0 & 134.0          & 232.0   & 67.0             \\ 
\hline
100                & 173.4 & 167.4         & 105.0 & 188.0            \\ 
\hline
150                & 181.2 & 164.2         & 297.2 & 256.2          \\ 
\hline
200                & 138.4  &  174.5        &  245.0     & 269.5                \\ 
\hline
250                & 148.0    &  169.4      &   195.6    & 168.0                \\ 
\hline
300                & 189.0    &  172.3              & 221.0       & 199.0               \\ 
\hline
350                & 177.1  & 174.0              & 189.0       & 155.0                \\ 
\hline
400                & 189.2       &  174.0              & 148.0      & 129.0                \\ 
\hline
450                & 168.0       &   191.0            & 205.7       & 178.0                \\ 
\hline
500                &  154.0     &  168.0      & 217.1       & 187.0               \\
\hline
\end{tabular}
\end{table}
}

\section{Conclusions}
We addressed the robust scheduling problem $R|r_j|C_{\max\limits}$ with interval release dates. The minimax regret criterion has been considered as the measure of robustness.  Our research includes the theoretical analysis, algorithms development and numerical experiments. Combining some theoretical features with the different greedy strategies allowed to develop the three efficient algorithms. Our approaches are based on both the simple makespan criterion and robust counterpart. Computational testing confirmed the complementariness of the developed algorithms and identified their running time limitations.\linebreak A substantial part of this study concerns the comparison of problem decomposition strategies. The results showed that the simplified and computational-effective strategies implemented in the PM an PR can provide competitive schedules as compared to the most complex approach used in the PRE. In fact, the PR ensures a fair trade off between the solution quality and running time. We have formulated two conditions when the makespan criterion allows us to obtain the optimal robust schedule.  Finally, the straightforward implementation is an unquestionable advantage of our algorithms.

Apart from the self-contained meaning of the investigated uncertain problem, it can be used as a tool for modeling and solving complex, fully deterministic optimization problems. Unlike the two-staged robust problem in \cite{liu2020parallel}, it can concern the joint scheduling and location problem (ScheLoc; e.g., \cite{lawrynowicz2019memetic}, \cite{kramer2021exact}, \cite{hessler2017discrete}).

\bibliographystyle{unsrt}  
\bibliography{manuscript}

\end{document}